\documentclass[12pt]{article}

\title{Non-parametric estimation in a  semimartingale regression
model.\\
Part 1. Oracle Inequalities.
\thanks{The paper is  supported by the RFFI-Grant 09-01-00172-a.}
}

\author{Konev Victor
\thanks{
Department of Applied Mathematics and Cybernetics,
 Tomsk State University,
Lenin str. 36,
 634050 Tomsk, Russia,
 e-mail: vvkonev@mail.tsu.ru }
 and
Pergamenshchikov Serguei\thanks{
 Laboratoire de Math\'ematiques Raphael Salem,
 Avenue de l'Universit\'e, BP. 12,
  Universit\'e de Rouen,
   F76801, Saint Etienne du Rouvray, Cedex France
and Department of Mathematics and Mechanics,Tomsk State University,
Lenin str. 36, 634041 Tomsk, Russia, e-mail:
Serge.Pergamenchtchikov@univ-rouen.fr } }

\usepackage{amssymb}
\usepackage{amsfonts}
\usepackage{amsmath}
\usepackage{amsthm}
\usepackage{enumerate}
\usepackage{multicol}

\newtheorem{theorem}{Theorem}[section]
\newtheorem{proposition}[theorem]{Proposition}
\newtheorem{lemma}[theorem]{Lemma}

\newtheorem{remark}{Remark}[section]
\newtheorem{corollary}[theorem]{Corollary}

\newcommand\cA{{\cal A}}

\newcommand\cH{{\cal H}}
\newcommand\cG{{\cal G}}

\newcommand\cL{{\cal L}}
\newcommand\cB{{\cal B}}

\newcommand\cD{{\cal D}}

\newcommand\cR{{\cal R}}

\newcommand\ve{\varepsilon}
\newcommand\ov{\overline}

\def\bbr{{\mathbb R}}

\def\text#1{\hbox{#1}}
\def\proof{{\noindent \bf Proof. }}
\def\endproof{\mbox{\ $\qed$}}
\def\card{\mbox{card}}

\def\E{{\bf E}}
\def\P{{\bf P}}
\def\C{{\bf C}}

\def\L{{\bf L}}

\newcommand{\wh}{\widehat}
\newcommand{\wt}{\widetilde}

\newcommand\Er{\mbox{Err}}

\def\Chi{{\bf 1}}
\def\d{\mathrm{d}}
\def\build #1_#2{\mathrel{\mathop{\kern 0pt #1}\limits_{#2}}}
\newcommand{\zs}[1]{{\mathchoice{#1}{#1}{\lower.25ex\hbox{$\scriptstyle#1$}}
{\lower0.25ex\hbox{$\scriptscriptstyle#1$}}}}
\numberwithin{equation}{section}

\begin{document}

\maketitle

\begin{abstract}

This paper considers the problem of estimating a periodic function
in a continuous time regression model with a general square integrable
semimartingale noise. A model
selection adaptive procedure is proposed. Sharp non-asymptotic
oracle inequalities have been derived.
\end{abstract}
\vspace*{5mm}
\noindent {\sl Keywords}: Non-asymptotic estimation; Non-parametric regression;
Model selection; Sharp oracle inequality; Semimartingale noise.

\vspace*{5mm}
\noindent {\sl AMS 2000 Subject Classifications}: Primary: 62G08; Secondary: 62G05

\bibliographystyle{plain}
\renewcommand{\columnseprule}{.1pt}

\newpage

\section{Introduction}\label{sec:In}

Consider a regression model in continuous time
 \begin{equation}\label{sec:In.1}
  \d y_t=S(t)\d t+\d \xi_t\,,
  \quad 0\le t\le n\,,
 \end{equation}
where $S$ is an unknown $1$-periodic $\bbr\to\bbr$ function, $S\in
\cL_\zs{2}[0,n]$; $(\xi_\zs{t})_\zs{t\ge 0}$ is a square
integrable unobservable semimartingale noise such that for any
function $f$ from $\cL_\zs{2}[0,n]$ the stochastic integral
\begin{equation}\label{sec:In.1-1}
I_\zs{n}(f)=\int^n_\zs{0}f_\zs{s}\d \xi_\zs{s}
 \end{equation} is
well defined with
\begin{equation}\label{sec:In.2}
\E I_\zs{n}(f)=0 \quad\mbox{and}\quad \E I^2_\zs{n}(f)\le
\sigma^*\,\int^n_\zs{0}\,f^2_\zs{s}\,\d s
\end{equation}
where $\sigma^*$ is some positive constant.

An important example of the disturbance $(\xi_\zs{t})_\zs{t\ge 0}$
is the following process

\begin{equation}\label{sec:In.3}
\xi_\zs{t}=\varrho_\zs{1}w_\zs{t}+\varrho_\zs{2}z_\zs{t}
\end{equation}
where $\varrho_\zs{1}$ and $\varrho_\zs{2}$ are unknown constants,
$|\varrho_\zs{1}|+|\varrho_\zs{2}|>0$, $(w_\zs{t})_\zs{t\ge 0}$ is a standard Brownian motion,
$(z_\zs{t})_\zs{t\ge 0}$ is a compound Poisson process  defined as
\begin{equation}\label{sec:In.4}
z_\zs{t}=\sum^{N_\zs{t}}_\zs{j=1}\,Y_\zs{j}
\end{equation}
where $(N_\zs{t})_\zs{t\ge 0}$ is a standard homogeneous Poisson process with
unknown  intensity $\lambda>0$ and
$(Y_\zs{j})_\zs{j\ge 1}$ is an i.i.d. sequence of random variables with
\begin{equation}\label{sec:In.5}
\E Y_\zs{j}=0
\quad\mbox{and}\quad
\E Y^2_\zs{j}=1\,.
\end{equation}
Let $(T)_\zs{k\ge 1}$ denote the arrival times of the process $(N_\zs{t})_\zs{t\ge 0}$,
that is,
\begin{equation}\label{sec:In.6}
T_\zs{k}=\inf\{t\ge 0\,:\,N_\zs{t}=k\}\,.
\end{equation}
As is shown in Lemma~\ref{Le.sec:A.2}, the condition
\eqref{sec:In.2} holds for the noise \eqref{sec:In.3} with
$\sigma^*=\varrho^2_\zs{1}+\varrho^2_\zs{2}\lambda$.

The problem is to estimate the unknown function $S$ in the model
\eqref{sec:In.1} on the basis of observations $(y_\zs{t})_\zs{0\le
t\le n}$.

This problem enables one to solve that of
functional statistics which is stated as follows.
 Let observations
$(x^{k})_\zs{0\le k\le n}$
be
a segment of a sequence of independent identically distributed
random processes
$x^{k}=(x^{k}_\zs{t})_\zs{0\le t\le 1}$
specified
on the interval $[0,1]$, which
obey the stochastic differential equations
\begin{equation}\label{sec:In.7}
\d x^{k}_\zs{t}=S(t)\d t+\d \xi^{k}_\zs{t}\,,\quad x^{k}_\zs{0}=x_\zs{0}\,,
\quad 0\le t\le 1\,,
\end{equation}
where $(\xi^{k})_\zs{1\le k\le n}$ is an i.i.d sequence of random
processes $\xi^{k}=(\xi^{k}_\zs{t})_\zs{0\le t\le 1}$ with the
same distribution as the process  \eqref{sec:In.3}. The problem is
to estimate the unknown function $f(t)\in\cL_\zs{2}[0,1]$ on the
basis of observations $(x^{k})_\zs{1\le k\le n}$.
 This model can be reduced to \eqref{sec:In.1}, \eqref{sec:In.3} in the following way.
Let $y=(y_\zs{t})_\zs{0\le t\ge n}$ denote
the process defined as :
$$
y_\zs{t}=
\left
\{
\begin{array}{ll}
x^{1}_\zs{t}\,,&\quad\mbox{if}\quad 0\le t\le 1\,;\\[2mm]
y_\zs{k-1}+x^{k}_\zs{t-k+1}-x_\zs{0}\,,&
\quad\mbox{if}\quad k-1\le t\le k\,,\quad 2\le k\le n\,.
\end{array}
\right.
$$
This process satisfies the stochastic
 differential equation
$$
\d y_\zs{t}=\wt{S}(t)\d t+\d \wt{\xi}_\zs{t}\,,
$$
where $\wt{S}(t)=S(\{t\})$ and
$$
\wt{\xi}_\zs{t}=
\left
\{
\begin{array}{ll}
\xi^{1}_\zs{t}\,,&\quad\mbox{if}\quad 0\le t\le 1\,;\\[2mm]
\wt{\xi}_\zs{k-1}+\xi^{k}_\zs{t-k+1}\,,&
\quad\mbox{if}\quad k-1\le t\le k\,,\quad 2\le k\le n\,;
\end{array}
\right.
$$
$\{t\}=t-[t]$ is the fractional part of number $t$.

In this paper we will consider the estimation problem for the
regression model \eqref{sec:In.1} in $\cL_\zs{2}[0,1]$ with
the quality of an estimate $\wh{S}$ being measured by the mean
integrated squared error (MISE)

\begin{equation}\label{sec:In.8}
\cR(\wh{S},S):=
\E_\zs{S}\,\|\wh{S}-S\|^2\,,
\end{equation}
where $\E_\zs{S}$ stands for the expectation with respect to the
distribution $\P_\zs{S}$ of the process \eqref{sec:In.1} given $S$;
$$
\|f\|^2:=\int^1_\zs{0}f^2(t) \d t\,.
$$

It is natural to treat this problem from the standpoint of the model selection approach.
The origin of this method goes back to early seventies with
the pioneering papers by Akaike \cite{Ak} and Mallows \cite{Ma} who proposed to introduce
penalizing in a log-likelihood type criterion.
 The further progress has been made by Barron, Birge and Massart \cite{BaBiMa}, \cite{Mas}
who developed a non-asymptotic model selection method which enabled one to derive
non-asymptotic oracle inequalities for
a gaussian non-parametric regression model with the  i.i.d. disturbance.
An oracle inequality yields the upper bound for the estimate risk via the minimal
risk corresponding to a chosen family of estimates.
 Galtchouk and Pergamenshchikov \cite{GaPe1} developed the Barron-Birge-Massart
technique treating the problem of estimating a non-parametric drift
function in a diffusion process from the standpoint of sequential analysis.
 Fourdrinier and Pergamenshchikov \cite{FoPe} extended the
 Barron-Birge-Massart method to the models with dependent observations and,
in contrast to all above-mentioned papers on the model selection method,
where the estimation procedures were based on the least squares estimates,
they proposed to use an arbitrary family of projective estimates in
an adaptive estimation procedure, and they discovered that one can employ
 the improved least square estimates to increase the estimation quality.
Konev and Pergamenshchikov \cite{KoPe2} applied this improved
model selection method to the non-parametric estimation problem of
a periodic function in a model with a coloured noise in continuous
time having
unknown spectral characteristics.
In all cited papers the non-asymptotic oracle inequalities
have been derived which
 enable one to establish the optimal convergence
rate for the minimax risks. Moreover, in the latter paper the
oracle inequalities have been found for the robust risks.

In addition to the optimal convergence rate, an important problem is that of the efficiency
of adaptive estimation procedures. In order to examine the efficiency property one has to obtain
 the oracle inequalities in which the principal term has the factor close to unity.

The first result in this direction is most likely due to Kneip \cite{Kn}
who obtained, for a gaussian regression model, the oracle inequality
with the factor close to unity at the principal term. The oracle
inequalities of this type were obtained as well in \cite{CaGo} and in
\cite{CaGoPiTs} for the inverse problems. It will be observed that the
derivation of oracle inequalities in all these papers rests upon
the fact that by applying the Fourier transformation one can
reduce the initial model to the statistical gaussian model with
independent observations. Such a transform is possible only for
gaussian models with independent homogeneous observations or for
the inhomogeneous ones with the known correlation characteristics.
This restriction significantly narrows the area of application of such
estimation procedures and rules out a broad class of models
including, in particular, widely used in econometrics
heteroscedastic regression models (see, for example, \cite{GoQu}).
For constructing adaptive procedures in the case of inhomogeneous
observations one needs to amend the approach to the estimation
problem. Galtchouk and Pergamenshchikov
 \cite{GaPe2}-\cite{GaPe3} have developed a new estimation method intended for the heteroscedastic
 regression models. The heart of this method is to combine the
 Barron-Birg\'e-Massart non-asymptotic penalization method \cite{BaBiMa} and the
 Pinsker weighted least square method
minimizing the asymptotic risk
  (see, for example, \cite{Nu}, \cite{Pi}). Combining of these approaches results
  in the significant improvement of the estimation quality
  (see numerical example in \cite{GaPe2}). As was shown in
\cite{GaPe3} and \cite{GaPe4}, the Galthouk-Pergamenshchikov
procedure is efficient with respect to the robust minimax risk,
i.e. the minimax risk with the additional supremum operation over
the whole family of addmissible model distributions. In the sequel
\cite{GaPe5}, \cite{GaPe6}, this approach has been applied to the
problem of a drift estimation in a  diffusion process. In this
paper we apply this procedure to the estimation of a regression
function $S$ in a semimartingale regression model \eqref{sec:In.1}. The rest
of the paper is organized as follows. In Section~\ref{sec:Mo} we
construct the model selection procedure on the basis of
weighted least squares estimates
 and state the main results in the form of oracle inequalities for the quadratic
 risks. Section~\ref{sec:Pr} gives the proofs of all theorems. In Appendix some
technical results are established.

\section{Model selection}\label{sec:Mo}

This Section gives the construction of a model selection procedure  for
estimating a function $S$ in \eqref{sec:In.1} on the basis of weighted least square estimates and states
the main results.

For estimating the unknown function $S$ in model \eqref{sec:In.1}, we
apply its Fourier expansion in the trigonometric basis
  $(\phi_j)_\zs{j\ge 1}$ in $\cL_2[0,1]$
defined as
\begin{equation}\label{sec:Mo.1}
\phi_1=1\,,\quad
\phi_\zs{j}(x)=\sqrt{2}\,Tr_\zs{j}(2\pi [j/2]x)\,,\ j\ge 2\,,
\end{equation}
where the function $Tr_\zs{j}(x)=\cos(x)$ for even $j$ and
$Tr_\zs{j}(x)=\sin(x)$ for odd $j$; $[x]$ denotes the integer part
of $x$. The corresponding Fourier coefficients
\begin{equation}\label{sec:Mo.2}
\theta_\zs{j}=(S,\phi_j)= \int^1_\zs{0}\,S(t)\,\phi_\zs{j}(t)\,\d
t
\end{equation}
can be estimated as
\begin{equation}\label{sec:Mo.3}
\wh{\theta}_\zs{j,n}= \frac{1}{n}\int^n_\zs{0}\,\phi_j(t)\,\d
y_\zs{t}\,.
\end{equation}
In view of \eqref{sec:In.1}, we obtain
\begin{equation}\label{sec:Mo.4}
\wh{\theta}_\zs{j,n}=\theta_\zs{j}+\frac{1}{\sqrt{n}}\xi_\zs{j,n}\,,
\quad
 \xi_\zs{j,n}=\frac{1}{\sqrt{n}}
I_\zs{n}(\phi_\zs{j})
\end{equation}
where $I_\zs{n}$ is given in \eqref{sec:In.1-1}.

For any sequence $x=(x_\zs{j})_\zs{j\ge 1}$, we set
\begin{equation}\label{sec:Mo.5}
|x|^2=\sum^\infty_\zs{j=1}x^2_\zs{j}
\quad\mbox{and}\quad
\#(x)=\sum^\infty_\zs{j=1}\,\Chi_\zs{\{|x_\zs{j}|>0\}}\,.
\end{equation}
Now we impose the additional conditions on the noise
$(\xi_\zs{t})_\zs{t\ge 0}$.\\[2mm]

\noindent $\C_\zs{1})$ {\it There exists some positive constant
$\sigma> 0$ such that the sequence
$$\varsigma_\zs{j,n}=\E
\xi^2_\zs{j,n}-\sigma\,,\quad j\ge 1\,,
$$
for any $n\ge 1$,
 satisfies the following
inequality
$$
c^*_\zs{1}(n)=
\sup_\zs{x\in\cH\,,\,\#(x)\le n}\,|B_\zs{1,n}(x)|\,
\,<\,\infty
$$
where $\cH=[-1,1]^{\infty}$ and
\begin{equation}\label{sec:Mo.5-1}
B_\zs{1,n}(x)=
\sum^\infty_\zs{j=1}\,x_\zs{j}\,\varsigma_\zs{j,n}\,.
\end{equation}
}
\vspace*{4mm}
\noindent $\C_\zs{2})$ {\it Assume that for all $n\ge 1$
$$
c^*_\zs{2}(n)= \sup_\zs{|x|\le 1\,,\#(x)\le n}\,
\E\,B^{2}_\zs{2,n}(x)\,<\,\infty
$$
where
\begin{equation}\label{sec:Mo.5-2}
B_\zs{2,n}(x)=
\sum^\infty_\zs{j=1}\,x_\zs{j}\,\wt{\xi}_\zs{j,n}
\quad\mbox{with}\quad
\wt{\xi}_\zs{j,n}=\xi^2_\zs{j,n}-\E \xi^2_\zs{j,n}\,.
\end{equation}
}
\vspace{2mm}

As is stated in Theorem~\ref{Th.sec:2.2},  Conditions $\C_\zs{1})$
and $\C_\zs{2})$ hold for the process \eqref{sec:In.3}.
 Further we
introduce a class of weighted least squares estimates for $S(t)$
defined as
\begin{equation}\label{sec:Mo.6}
\wh{S}_\zs{\gamma}=\sum^{\infty}_\zs{j=1}\gamma(j)\wh{\theta}_\zs{j,n}\phi_\zs{j}\,,
\end{equation}
where $\gamma=(\gamma(j))_\zs{j\ge 1}$ is
a sequence of weight coefficients such that
\begin{equation}\label{sec:Mo.6-1}
0\le \gamma(j)\le 1 \quad\mbox{and}\quad 0<\#(\gamma)\le n\,.
\end{equation}
 Let $\Gamma$ denote a finite set of weight sequences $\gamma=(\gamma(j))_\zs{j\ge 1}$ with
 these properties, $\nu=\card(\Gamma)$ be its cardinal number and
\begin{equation}\label{sec:Mo.7}
\mu=\max_\zs{\gamma\in\Gamma}\,\#(\gamma)\,.
\end{equation}
\noindent The model selection procedure for the unknown function
$S$ in \eqref{sec:In.1} will be constructed on the basis of
 estimates $(\wh{S}_\zs{\gamma})_\zs{\gamma\in\Gamma}$.
The choice of a specific set
of weight sequences
 $\Gamma$ will be discussed at the end
of this section. In order to find a proper weight sequence
$\gamma$ in the set $\Gamma$ one needs  to specify a cost function.
When choosing an appropriate cost function one can use the
following argument. The empirical squared error
$$
\Er_\zs{n}(\gamma)=\|\wh{S}_\zs{\gamma}-S\|^2
$$
can be written as
\begin{equation}\label{sec:Mo.8}
\Er_\zs{n}(\gamma)\,=\,
\sum^\infty_\zs{j=1}\,\gamma^2(j)\wh{\theta}^2_\zs{j,n}\,-
2\,\sum^\infty_\zs{j=1}\,\gamma(j)\wh{\theta}_\zs{j,n}\,\theta_\zs{j}\,+\,
\sum^\infty_\zs{j=1}\theta^2_\zs{j}\,.
\end{equation}
Since the Fourier coefficients $(\theta_\zs{j})_\zs{j\ge 1}$ are
unknown, the weight coefficients $(\gamma_\zs{j})_\zs{j\ge 1}$ can not be
determined by minimizing this quantity. To circumvent this
difficulty one needs to replace  the terms
$\wh{\theta}_\zs{j,n}\,\theta_\zs{j}$ by some their estimators
$\wt{\theta}_\zs{j,n}$. We set
\begin{equation}\label{sec:Mo.9}
\wt{\theta}_\zs{j,n}=
\wh{\theta}^2_\zs{j,n}-\frac{\wh{\sigma}_\zs{n}}{n}
\end{equation}
where $\wh{\sigma}_\zs{n}$ is an estimator for the quantity $\sigma$ in
condition $\C_\zs{1})$.

For this change in the empirical squared error, one has to pay
some penalty. Thus, one comes to the cost function of the form
\begin{equation}\label{sec:Mo.10}
J_\zs{n}(\gamma)\,=\,\sum^\infty_\zs{j=1}\,\gamma^2(j)\wh{\theta}^2_\zs{j,n}\,-
2\,\sum^\infty_\zs{j=1}\,\gamma(j)\,\wt{\theta}_\zs{j,n}\,
+\,\rho\,\wh{P}_\zs{n}(\gamma)
\end{equation}
where $\rho$ is some positive constant,
$\wh{P}(\gamma)$ is the penalty term defined as
\begin{equation}\label{sec:Mo.11}
\wh{P}_\zs{n}(\gamma)=\frac{\wh{\sigma}_\zs{n}\,|\gamma|^2}{n}
\,.
\end{equation}

\noindent In the case when the value of $\sigma$ in $\C_\zs{1})$ is known, one can put
$\wh{\sigma}_\zs{n}=\sigma$ and
\begin{equation}\label{sec:Mo.12}
P_\zs{n}(\gamma)=\frac{\sigma\,|\gamma|^2}{n}\,.
\end{equation}
Substituting the weight coefficients, minimizing the cost function, that is
\begin{equation}\label{sec:Mo.13}
\wh{\gamma}=\mbox{argmin}_\zs{\gamma\in\Gamma}\,J_n(\gamma)\,,
\end{equation}
in \eqref{sec:Mo.6} leads to the model selection procedure
\begin{equation}\label{sec:Mo.14}
\wh{S}_\zs{*}=\wh{S}_\zs{\wh{\gamma}}\,.
\end{equation}
It will be noted that
$\wh{\gamma}$ exists, since
 $\Gamma$ is a finite set. If the
minimizing sequence in \eqref{sec:Mo.13} $\wh{\gamma}$ is not
unique, one can take any minimizer.

\begin{theorem}\label{Th.sec:2.1}
Assume that the conditions $\C_\zs{1})$ and $\C_\zs{2})$ hold with
$\sigma>0$. Then for any $n\ge 1$ and $0<\rho<1/3$, the estimator
\eqref{sec:Mo.14} satisfies the oracle inequality
\begin{align}\label{sec:Mo.15}
\cR(\wh{S}_\zs{*},S)\,\le\, \frac{1+3\rho-2\rho^2}{1-3\rho}
\min_\zs{\gamma\in\Gamma} \cR(\wh{S}_\zs{\gamma},S)
+\frac{1}{n}\,\cB^{*}_\zs{n}(\rho)
\end{align}
where the risk $\cR(\cdot,S)$ is defined in \eqref{sec:In.8},
$$
\cB^{*}_\zs{n}(\rho)=\Psi_\zs{n}(\rho)+
\frac{6\mu\,
\E_\zs{S}|\wh{\sigma}_\zs{n}-\sigma|}{1-3\rho}
$$
and
\begin{equation}\label{sec:Mo.15-1}
\Psi_\zs{n}(\rho)= \frac{ 2\sigma\sigma^*\nu +
4\sigma c^*_\zs{1}(n) + 2\nu c^{*}_\zs{2}(n)}{\sigma\rho(1-3\rho)}
\,.
\end{equation}
\end{theorem}
Now we check conditions $\C_\zs{1})$ and $\C_\zs{2})$ for the model
\eqref{sec:In.1} with the noise  \eqref{sec:In.3} to arrive at the following result.

\vspace{2mm}

\begin{theorem}\label{Th.sec:2.2}
Suppose that  the coefficients $\varrho_\zs{1}$ and
$\varrho_\zs{2}$ in model \eqref{sec:In.1}, \eqref{sec:In.3},
 are such that
$\varrho^2_\zs{1}+\varrho^2_\zs{2}>0$ and $\E Y^4_\zs{j}<\infty$.
Then the estimation procedure \eqref{sec:Mo.14}, for any $n\ge 1$
and $0<\rho\le 1/3$,
 satisfies the oracle inequality
\eqref{sec:Mo.15} with
$$
\sigma=
\sigma^*=\varrho^2_\zs{1}+\lambda\varrho^2_\zs{2}
\,,\quad
c^*_\zs{1}(n)=0\,,
$$
and
$$
\sup_\zs{n\ge 1}\, c^*_\zs{2}(n)
\le\,4
\sigma^{*}
\left(\sigma^* + \varrho^2_\zs{2}\E
Y^4_\zs{1} \right)
\,.
$$
\end{theorem}
\noindent The proofs of Theorems~\ref{Th.sec:2.1},
\ref{Th.sec:2.2} are given in Section~\ref{sec:Pr}.

\begin{corollary}\label{Co.sec:2.1}
Let the conditions of Theorem~\ref{Th.sec:2.1} hold and the
quantity
$\sigma$ in $\C_\zs{1})$ be known. Then, for any $n\ge 1$
and $0<\rho<1/3$, the estimator \eqref{sec:Mo.14} satisfies the
oracle inequality
$$
\cR(\wh{S}_\zs{*},S)\,\le\, \frac{1+3\rho-2\rho^2}{1-3\rho}
\min_\zs{\gamma\in\Gamma} \cR(\wh{S}_\zs{\gamma},S)
+\frac{1}{n}\,\Psi_\zs{n}(\rho)\,,
$$
where $\Psi_\zs{n}(\rho)$ is given in \eqref{sec:Mo.15-1}.
\end{corollary}

\subsection{Estimation of $\sigma$}\label{subsec:Si}

Now we consider the case of unknown quantity $\sigma$ in the condition
$\C_\zs{1})$. One can estimate $\sigma$ as

\begin{equation}\label{sec:Si.1}
\wh{\sigma}_\zs{n}=\sum^n_\zs{j=l}\,\wh{\theta}^2_\zs{j,n}
\quad\mbox{with}\quad l=[\sqrt{n}]+1\,.
\end{equation}
\vspace{2mm}
\begin{proposition}\label{Pr.sec:Si.1}
Suppose that the conditions of Theorem~\ref{Th.sec:2.1} hold and
the unknown function $S(t)$ is continuously differentiable
for $0\le t<1$ such that
\begin{equation}\label{sec:Si.2}
 |\dot{S}|_\zs{1}=
\int^1_\zs{0}|\dot{S}(t)|\d t < +\infty\,.
\end{equation}
Then, for any $n\ge 1$,
\begin{equation}\label{sec:Si.3}
\E_\zs{S}|\wh{\sigma}_\zs{n}-\sigma|\le
\frac{\kappa_\zs{n}(S)}{\sqrt{n}}
\end{equation}
where
$$
\kappa_\zs{n}(S)= 4|\dot{S}|^2_\zs{1}+\sigma+\sqrt{c^*_\zs{2}(n)}+
\frac{4|\dot{S}|_\zs{1}\sqrt{\sigma^*}}{n^{1/4}} +
\frac{c^*_\zs{1}(n)}{n^{1/2}}\,.
$$
\end{proposition}
The proof of Proposition~\ref{Pr.sec:Si.1} is given in Section~\ref{sec:Pr}.
Theorem~\ref{Th.sec:2.1} and Proposition~\ref{Pr.sec:Si.1}
imply the following result.
\begin{theorem}\label{Th.sec:Si.1}
Suppose that the conditions of Theorem~\ref{Th.sec:2.1} hold and
$S$ satisfies the conditions of Proposition~\ref{Pr.sec:Si.1}.
Then, for any $n\ge 1$ and
$0<\rho<1/3$, the estimate \eqref{sec:Mo.14} satisfies the oracle
inequality
\begin{align}\label{sec:Si.4}
\cR(\wh{S}_\zs{*},S)\,\le\, \frac{1+3\rho-2\rho^2}{1-3\rho}
\min_\zs{\gamma\in\Gamma} \cR(\wh{S}_\zs{\gamma},S)
+\frac{1}{n}\,\cD_\zs{n}(\rho)\,,
\end{align}
where
$$
\cD_\zs{n}(\rho) =2\,\Psi_\zs{n}(\rho)+ \frac{2\rho(1-\rho)\mu
\kappa_\zs{n}(S)}{(1-3\rho)\,\sqrt{n}}
 \,.
$$
\end{theorem}

\subsection{Specification of weights in the selection procedure \eqref{sec:Mo.14}}\label{subsec:Ga}

Now we will specify the weight coefficients $(\gamma(j))_\zs{j\ge
1}$ in a way proposed in \cite{GaPe2}
 for a heteroscedastic discrete time regression
model. Consider a numerical grid of the form
$$
\cA_\zs{n}=\{1,\ldots,k^*\}\times\{t_1,\ldots,t_m\}\,,
$$
where  $t_i=i\ve$ and $m=[1/\ve^2]$. We assume that both
parameters $k^*\ge 1$ and $0<\ve\le 1$ are functions of $n$, i.e.
$ k^*=k^*(n)$ and $\ve=\ve(n)$, such that
\begin{equation}\label{sec:Ga.1}
\left\{
\begin{array}{ll}
&\lim_\zs{n\to\infty}\,k^*(n)=+\infty\,,
\quad
\lim_\zs{n\to\infty}\,\dfrac{k^*(n)}{\ln n}=0\,,\\[6mm]
&\lim_\zs{n\to\infty}\ve(n)=0
\quad\mbox{and}\quad
\lim_\zs{n\to\infty}\,n^{\delta}\ve(n)\,=+\infty
\end{array}
\right.
\end{equation}
for any $\delta>0$. One can take, for example,
$$
\ve(n)=\frac{1}{\ln (n+1)}
\quad\mbox{and}\quad
k^*(n)=\sqrt{\ln (n+1)}
$$
for  $n\ge 1$.

For each $\alpha=(\beta,t)\in\cA_\zs{n}$, we introduce the weight
sequence
$$
\gamma_\zs{\alpha}=(\gamma_\zs{\alpha}(j))_\zs{j\ge 1}
$$
given as
\begin{equation}\label{sec:Ga.2}
\gamma_\zs{\alpha}(j)=\Chi_\zs{\{1\le j\le j_\zs{0}\}}+
\left(1-(j/\omega_\alpha)^\beta\right)\, \Chi_\zs{\{ j_\zs{0}<j\le
\omega_\alpha\}}
\end{equation}
where $j_\zs{0}=j_\zs{0}(\alpha)=\left[\omega_\zs{\alpha}/\ln n\right]$,
$$
\omega_\zs{\alpha}=(\tau_\zs{\beta}\,t\,n)^{1/(2\beta+1)}
\quad\mbox{and}\quad
\tau_\zs{\beta}=\frac{(\beta+1)(2\beta+1)}{\pi^{2\beta}\beta}\,.
$$
We set
\begin{equation}\label{sec:Ga.3}
\Gamma\,=\,\{\gamma_\zs{\alpha}\,,\,\alpha\in\cA_\zs{n}\}\,.
\end{equation}
It will be noted that in this case $\nu=k^* m$.

\begin{remark}\label{Re.sec:Ga.1}
It will be observed that the specific form of weights
\eqref{sec:Ga.2} was proposed by Pinsker \cite{Pi} for the
filtration problem with known smoothness of regression function
observed with an additive gaussian white noise in the continuous
time. Nussbaum \cite{Nu} used these  weights for the
gaussian regression estimation problem in discrete time.

The minimal mean square risk,  called the Pinsker
constant, is provided by the weight least squares estimate with
the weights where the index $\alpha$ depends on the smoothness
order of the function $S$. In this case the smoothness order is
unknown and, instead of one estimate, one has to use a whole family
of estimates containing in particular the optimal one.

The problem is to study the properties of the whole class of
estimates. Below we derive an oracle inequality for this class
which yields the best mean square risk up to a multiplicative and
additive constants
provided that the
 the smoothness of the unknown function $S$
is not available. Moreover, it will be shown that the
multiplicative constant tends to unity and the additive one vanishes
as $n\to\infty$ with the rate higher than any minimax rate.
\end{remark}
\noindent In view of the assumptions \eqref{sec:Ga.1}, for any $\delta>0$,
one has
$$
\lim_\zs{n\to\infty}\,\frac{\nu}{n^{\delta}}=0\,.
$$
Moreover, by \eqref{sec:Ga.2} for any $\alpha\in\cA_\zs{n}$
\begin{align*}
\sum^\infty_\zs{j=1}\,\Chi_\zs{\{\gamma_\zs{\alpha}(j)>0\}}
  \le \omega_\zs{\alpha}\,.
\end{align*}
Therefore, taking into account that $A_\zs{\beta}\le A_1<1$ for
$\beta\ge 1$, we get
$$
\mu=\mu_\zs{n}\le(n/\ve )^{1/3}\,.
$$
Therefore, for any $\delta>0$,
$$
\lim_\zs{n\to\infty}\frac{\mu_\zs{n}}{n^{1/3+\delta}}=0\,.
$$
Applying this limiting relation to the analysis of the asymptotic
behavior of the additive term $\cD_\zs{n}(\rho)$ in
\eqref{sec:Si.4} one comes to the following result.

\begin{theorem}\label{Th.sec:La.1}
Suppose that the conditions of Theorem~\ref{Th.sec:2.1} hold and
$\dot{S}\in \cL_\zs{1}[0,1]$.  Then, for any $n\ge 1$ and
$0<\rho<1/3$, the estimate \eqref{sec:Mo.14} with the weight
coefficients \eqref{sec:Ga.3} satisfies the oracle inequality
\eqref{sec:Si.4} with the additive term $\cD_\zs{n}(\rho)$
obeying, for any $\delta>0$, the following limiting relation
$$
\lim_\zs{n\to\infty}\,\frac{\cD_\zs{n}(\rho)}{n^\delta}=0
 \,.
$$
\end{theorem}

\vspace*{5mm}
\section{Proofs}\label{sec:Pr}
\subsection{Proof of Theorem~\ref{Th.sec:2.1}}\label{subsec:Pr.1}

Substituting \eqref{sec:Mo.10} in \eqref{sec:Mo.8} yields for any $\gamma\in\Gamma$

\begin{equation}\label{sec:Pr.1}
\Er_\zs{n}(\gamma)=J_\zs{n}(\gamma)+
2\sum^\infty_\zs{j=1}\gamma(j)\theta'_\zs{j,n}+
\|S\|^2-\rho\,
\wh{P}_\zs{n}(\gamma)\,,
\end{equation}
where
$$
\theta'_\zs{j,n}=
\wt{\theta}_\zs{j,n}\,-\,\theta_\zs{j}\wh{\theta}_\zs{j,n}\\
=\frac{1}{\sqrt{n}}\theta_\zs{j}\xi_\zs{j,n}+
\frac{1}{n}\wt{\xi}_\zs{j,n}
+\frac{1}{n}\varsigma_\zs{j,n}
+\frac{\sigma-\wh{\sigma}_\zs{n}}{n}
$$
and the sequences $(\varsigma_\zs{j,n})_\zs{j\ge 1}$
and $(\wt{\xi}_\zs{j,n})_\zs{j\ge 1}$ are defined in conditions $\C_\zs{1})$
and $\C_\zs{2})$.
Denoting
\begin{equation}\label{sec:Pr.2}
L(\gamma)=\sum^\infty_\zs{j=1}\,\gamma(j)\,,\quad
M(\gamma)=
\frac{1}{\sqrt{n}}
\sum^\infty_\zs{j=1}\gamma(j)\theta_\zs{j}\xi_\zs{j,n}\,,
\end{equation}
and taking into account the definition of the "true" penalty term in \eqref{sec:Mo.12},
we rewrite \eqref{sec:Pr.1} as
\begin{align}\nonumber
\Er_\zs{n}(\gamma)&=J_\zs{n}(\gamma)
+2\frac{\sigma-\wh{\sigma}_\zs{n}}{n}\,L(\gamma)+
2M(\gamma) +\frac{2}{n} B_\zs{1,n}(\gamma)
\\[2mm] \label{sec:Pr.3}
&+2\sqrt{P_\zs{n}(\gamma)}
\frac{B_\zs{2,n}(e(\gamma))}{\sqrt{\sigma n}}
+\|S\|^2 -\rho
\wh{P}_\zs{n}(\gamma)
\end{align}
where $e(\gamma)=\gamma/|\gamma|$, the functions $B_\zs{1,n}$ and $B_\zs{2,n}$
are defined in \eqref{sec:Mo.5-1} and \eqref{sec:Mo.5-2}.

Let $\gamma_\zs{0}=(\gamma_\zs{0}(j))_\zs{j\ge 1}$ be a fixed
sequence in $\Gamma$ and $\wh{\gamma}$ be as in \eqref{sec:Mo.13}.
Substituting $\gamma_\zs{0}$ and $\wh{\gamma}$ in the equation
\eqref{sec:Pr.3}, we consider the difference
\begin{align*}
\Er_\zs{n}(\wh{\gamma})-\Er_\zs{n}(\gamma_\zs{0})&=J(\wh{\gamma})-J(\gamma_\zs{0})
+
2\frac{\sigma-\wh{\sigma}}{n}\,L(\wh{x})
+\frac{2}{n}
B_\zs{1,n}(\wh{x})
+
2M(\wh{x})
\\
&+ 2\sqrt{P_\zs{n}(\wh{\gamma})}
\frac{B_\zs{2,n}(\wh{e})}{\sqrt{\sigma n}}
-2\sqrt{P_\zs{n}(\gamma_\zs{0})} \frac{B_\zs{2,n}(e_\zs{0})}{\sqrt{\sigma n}}
\\
&
-\rho \wh{P}_\zs{n}(\wh{\gamma})
+\rho \wh{P}_\zs{n}(\gamma_\zs{0})
\end{align*}
where
$\wh{x}=\wh{\gamma}-\gamma_\zs{0}$, $\wh{e}=e(\wh{\gamma})$ and $e_\zs{0}=e(\gamma_\zs{0})$.
 Note that
by \eqref{sec:Mo.7}
$$
|L(\wh{x})|\le |L(\wh{\gamma})|+|L(\gamma)|\le 2\mu\,.
$$
Therefore, by making
use of the condition $\C_\zs{1})$
 and taking into account that the cost
function $J$ attains its minimum at $\wh{\gamma}$, one comes to
the inequality
\begin{align}\nonumber
\Er_\zs{n}(\wh{\gamma})-\Er_\zs{n}(\gamma_\zs{0})&\le
4\frac{|\wh{\sigma}-\sigma|}{n}\,\mu\,
+
\frac{2c^*_\zs{1}(n)}{n}
+
2M(\wh{x})
\\ \nonumber
&+2\sqrt{P_\zs{n}(\wh{\gamma})} \frac{B_\zs{2,n}(\wh{e})}{\sqrt{\sigma n}}
-\rho \wh{P}_\zs{n}(\wh{\gamma})
\\
\label{sec:Pr.5}
& +\rho \wh{P}_\zs{n}(\gamma_\zs{0})
-2\sqrt{P_\zs{n}(\gamma_\zs{0})}
\frac{B_\zs{2,n}(e_\zs{0})}{\sqrt{\sigma n}}\,.
\end{align}
Applying the elementary inequality
\begin{equation}\label{sec:Pr.6}
2|ab|\le \ve a^2+\ve^{-1}b^2
\end{equation}
with $\ve=\rho$ implies
the estimate
$$
2 \sqrt{P_\zs{n}(\gamma)}
\frac{|B_\zs{2,n}(e(\gamma))|}{\sqrt{\sigma n}}\le \rho P_\zs{n}(\gamma)
+\frac{B^2_\zs{2,n}(e(\gamma))}{n\sigma \rho}\,.
$$
We recall that $0<\rho<1$. Therefore, from here and \eqref{sec:Pr.5}, it follows
that
\begin{align*}
\Er_\zs{n}(\wh{\gamma})&\le\Er_\zs{n}(\gamma_\zs{0})+
2M(\wh{x})
+
\frac{2 B^*_\zs{2,n}}{n\sigma \rho }
+\frac{2c^*_\zs{1}(n)}{n}\\[2mm]
&
+
\frac{1}{n}|\wh{\sigma}-\sigma|
\left(
|\wh{\gamma}|^2+|\gamma_\zs{0}|^2+4\mu
\right)
+2\rho P_\zs{n}(\gamma_\zs{0})
\end{align*}
where $B^{*}_\zs{2,n}=\sup_\zs{\gamma\in\Gamma} B^{2}_\zs{2,n}(e(\gamma))$.
In view of \eqref{sec:Mo.7}, one has
$$
\sup_\zs{\gamma\in\Gamma}|\gamma|^2\le \mu\,.
$$
Thus, one gets
\begin{align}\nonumber
\Er_\zs{n}(\wh{\gamma})&\le\Er_\zs{n}(\gamma_\zs{0})+ 2M(\wh{x}) +
\frac{2 B^*_\zs{2,n}}{n\sigma \rho } +\frac{2c^*_\zs{1}(n)}{n}\\[2mm]
\label{sec:Pr.7}
& + \frac{6\mu}{n}|\wh{\sigma}_\zs{n}-\sigma| +2\rho
P_\zs{n}(\gamma_\zs{0})\,.
\end{align}
In view of Condition $\C_\zs{2})$, one has
\begin{align}\label{sec:Pr.8}
\E_\zs{S}\,B^{*}_\zs{2,n}\,
\le\,\sum_\zs{\gamma\in\Gamma}\,\E_\zs{S}\,B^{2}_\zs{2,n}(e(\gamma))
\,\le \,\nu\,c^*_\zs{2}(n)
\end{align}
where $\nu=\card(\Gamma)$.

Now we examine the first term in the right-hand side of
\eqref{sec:Pr.5}. Substituting \eqref{sec:Mo.4} in
\eqref{sec:Pr.2} and taking into account \eqref{sec:In.2}, one
obtains that for any non-random sequence $x=(x(j))_\zs{j\ge 1}$
with $\#(x)<\infty$
\begin{equation}\label{sec:Pr.9}
\E_\zs{S} M^2(x)
\le\,
\sigma^*
\,\frac{1}{n}\,
\sum^\infty_\zs{j=1}\,x^2(j)
\theta^2_\zs{j}\,
=
\,\sigma^*\,
\frac{1}{n}\,
\|S_\zs{x}\|^2
\end{equation}
where $S_\zs{x}=\sum^\infty_\zs{j=1}x(j)\theta_\zs{j}\phi_\zs{j}$.
Let denote
$$
Z^*=\sup_\zs{x\in \Gamma_\zs{1}}
\frac{n M^2(x)}{\|S_\zs{x}\|^2}
$$
where $\Gamma_\zs{1}=\Gamma-\gamma_\zs{0}$.
In view of \eqref{sec:Pr.9}, this quantity can be estimated as
\begin{equation}\label{sec:Pr.10}
\E_\zs{S}\,Z^*\,\le\,
\sum_\zs{x\in \Gamma_\zs{1}}
\frac{n\E_\zs{S}\,M^2(x)}{\|S_\zs{x}\|^2}
\le
\sum_\zs{x\in \Gamma_\zs{1}}\,\sigma^*=
 \sigma^*\nu\,.
\end{equation}
Further, by making use of the inequality \eqref{sec:Pr.6} with
$\varepsilon=\rho\|S_\zs{x}\|$, one gets
\begin{equation}\label{sec:Pr.11}
2|M(x)|\le\rho \|S_\zs{x}\|^2+
\frac{Z^*}{n\rho}\,.
\end{equation}
Note that, for any $x\in \Gamma_\zs{1}$,
\begin{equation}\label{sec:Pr.12}
\|S_\zs{x}\|^2-
\|\wh{S}_\zs{x}\|^2=
\sum^\infty_\zs{j=1}\,x^2(j)
(\theta^2_\zs{j}-\wh{\theta}^2_\zs{j,n})\le -2 M_\zs{1}(x)
\end{equation}
where
$$
M_\zs{1}(x)=
\frac{1}{\sqrt{n}}\sum^\infty_\zs{j=1}x^2(j)\theta_\zs{j}\xi_\zs{j,n}\,.
$$
Since $|x(j)|\le 1$
 for any $x\in\Gamma_\zs{1}$, one gets
$$
\E_\zs{S} M^2_\zs{1}(x)\le \sigma^*
\frac{\|S_\zs{x}\|^2}{n}\,.
$$
Denoting
$$
Z^*_\zs{1}=\sup_\zs{x\in \Gamma_\zs{1}}
\frac{nM^2_\zs{1}(x)}{\|S_\zs{x}\|^2}\,,
$$
one has
\begin{equation}\label{sec:Pr.13}
\E_\zs{S}\,Z^*_\zs{1}\,\le \sigma^*\nu\,.
\end{equation}
By the same argument as in \eqref{sec:Pr.11},
 one derives
 $$
2|M_\zs{1}(x)|\le\rho \|S_\zs{x}\|^2
+
\frac{Z^*_\zs{1}}{n\rho}\,.
$$
From here and \eqref{sec:Pr.12}, one finds the upper bound for
$\|S_\zs{x}\|$, i.e.
\begin{align}\label{sec:Pr.14}
\|S_\zs{x}\|^2&\le\,\frac{\|\wh{S}_\zs{x}\|^2}{1-\rho}
+\,
\frac{Z^*_\zs{1}}{n\rho(1-\rho)}\,.
\end{align}
Using this bound in \eqref{sec:Pr.11} gives
$$
2M(x)\le\,\frac{\rho \|\wh{S}_\zs{x}\|^2}{1-\rho }
+
\frac{Z^*+Z^*_\zs{1}}{n\rho(1-\rho)}\,.
$$
Setting
$x=\wh{x}$
in this inequality
 and taking into account that
$$
\|\wh{S}_\zs{\wh{x}}\|^2=\|\wh{S}_\zs{\wh{\gamma}}-\wh{S}_\zs{\gamma_\zs{0}}\|^2
\le
2
(\Er_n(\wh{\gamma})+\Er_n(\gamma_\zs{0}))\,,
$$
we obtain
$$
2M(\wh{x})
\le
\frac{2\rho(\Er_n(\wh{\gamma})+\Er_n(\gamma_\zs{0}))}{1-\rho}
+
\frac{Z^*+Z^*_\zs{1}}{n\rho(1-\rho)}\,.
$$
From here and \eqref{sec:Pr.7}, it follows that
\begin{align*}
\Er_\zs{n}(\wh{\gamma})&
\le \frac{1+\rho }{1-3\rho }
\Er_\zs{n}(\gamma_\zs{0})+
\frac{2(1-\rho) }{n(1-3\rho)}
\left(
\frac{ B^{*}_\zs{2,n}}{\sigma \rho}
+c^*_\zs{1}(n)
+
3\mu|\wh{\sigma}-\sigma|
\right)
\\[2mm]
&+
\frac{Z^*+Z^*_\zs{1}}{n\rho(1-3\rho )}
+\frac{2\rho(1-\rho)}{1-3\rho}P_\zs{n}(\gamma_\zs{0})\,,
\end{align*}
Taking the expectation yields
\begin{align*}
\cR(\wh{S}_\zs{*},S)
&
\le \frac{1+\rho }{1-3\rho }
\cR(\wh{S}_\zs{\gamma_\zs{0}},S)
+
\frac{2(1-\rho) }{n(1-3\rho)}
\left(
\frac{\nu c^*_\zs{2}(n)}{\sigma \rho}
+c^*_\zs{1}(n)
+
3\mu\,\E_\zs{S}\,
|\wh{\sigma}-\sigma|
\right)
\\[2mm]
&+
\frac{2\sigma^*\nu}{n\rho(1-3\rho )}
+\frac{2\rho(1-\rho)}{1-3\rho}P_\zs{n}(\gamma_\zs{0})\,.
\end{align*}
Using the upper bound for $P_\zs{n}(\gamma_\zs{0})$ in
Lemma~\ref{Le.sec:A.1}, one obtains
\begin{align*}
\cR(\wh{S}_\zs{*},S)\,\le\,
\frac{1+3\rho-2\rho^2}{1-3\rho}
\cR(\wh{S}_\zs{\gamma_\zs{0}},S)
+\frac{1}{n}\,\cB^{*}_\zs{n}(\rho)\,,
\end{align*}
where $\cB^{*}_\zs{n}(\rho)$ is defined in \eqref{sec:Mo.15}.

Since this inequality holds for each $\gamma_\zs{0}\in\Gamma$,
this completes the proof of
 Theorem~\ref{Th.sec:2.1}.
\endproof

\vspace{2mm}

\subsection{Proof of Theorem~\ref{Th.sec:2.2}}\label{subsec:Pr.2}
We have to verify Conditions $\C_\zs{1})$ and $\C_\zs{2})$ for the
process \eqref{sec:In.3}.

Condition $\C_\zs{1})$ holds with
$c^*_\zs{1}(n)=0$. This follows from Lemma~\ref{Le.sec:A.2} if one puts $f=g=\phi_\zs{j}$, $j\ge 1$.
Now we check Condition $\C_\zs{2})$.
By the Ito formula and Lemma~\ref{Le.sec:A.2}, one gets
$$
\d I^2_\zs{t}(f)=2I_\zs{t-}(f)\d I_\zs{t}(f)
+\varrho^2_\zs{1}f^2(t)\d t
+\varrho^2_\zs{2}\,
\sum_\zs{0\le s\le t}f^2(s)(\Delta z_\zs{s})^2
$$
and
$$
\E\,I^2_\zs{t}(f)=
\sigma^{*}
\int^{t}_\zs{0}
\,f^2(t)\d t\,.
$$
Therefore, putting
$$
\wt{I}_\zs{t}(f)=I^2_\zs{t}(f)-\E I^2_\zs{t}(f)\,,
$$
we obtain
$$
\d \wt{I}_\zs{t}(f)= \varrho^2_\zs{2}\,f^{2}(t)\,\d m_\zs{t}+
2I_\zs{t-}(f)f(t)\d \xi_\zs{t}\,, \quad \wt{I}_\zs{0}(f)=0
$$
 and
\begin{equation}\label{sec:Pr.16}
m_\zs{t} =  \sum_\zs{0\le s\le t}\, (\Delta z_\zs{s})^{2}
-\,\lambda\,t\,.
\end{equation}
Now we set
$$
\ov{I}_\zs{t}(x)=\sum^\infty_\zs{j=1}x_\zs{j}\wt{I}_\zs{t}(\phi_\zs{j})
$$
where $x=(x_\zs{j})_\zs{j\ge 1}$ with $\#(x)\le n$ and $|x|\le 1$.
This process obeys the equation
$$
\d \ov{I}_\zs{t}(x)= \varrho^2_\zs{2}\,\Phi_\zs{t}\,\d m_\zs{t} +
2\zeta_\zs{t-}(x)\d \xi_\zs{t} \,,\quad \ov{I}_\zs{0}(x)=0\,,
$$
where
$$
\Phi_\zs{t}(x)=
\sum_\zs{j\ge 1}x_\zs{j}\,\phi^2_\zs{j}(t)
\quad\mbox{and}\quad
\zeta_\zs{t}(x)
=
\sum_\zs{j\ge 1}x_\zs{j}I_\zs{t}(\phi_\zs{j})\phi_\zs{j}(t)
\,.
$$
Now we show that
\begin{equation}\label{sec:Pr.16-1}
\E\,\int^n_\zs{0}\,\ov{I}_\zs{t-}(x)\d \ov{I}_\zs{t}(x)=0\,.
\end{equation}
Indeed, note that
\begin{align*}
\int^n_\zs{0}\,\ov{I}_\zs{t-}(x)\d \ov{I}_\zs{t}(x)
&=
\varrho^{2}_\zs{2}
\sum_\zs{j\ge 1}\,x_\zs{j}
\int^{n}_\zs{0}\,\wt{I}_\zs{t-}(\phi_\zs{j})\Phi_\zs{t}(x)\d m_\zs{t}\\
&+
2
\sum_\zs{j\ge 1}\,x_\zs{j}
\int^{n}_\zs{0}\,\wt{I}_\zs{t-}(\phi_\zs{j})\zeta_\zs{t-}(x)\d\xi_\zs{t}
\,.
\end{align*}
Therefore, Lemma~\ref{Le.sec:A.4} directly implies
\begin{align*}
\E\,\int^{n}_\zs{0}\,\wt{I}_\zs{t-}(\phi_\zs{j})\Phi_\zs{t}(x)\d m_\zs{t}
&=
\sum_\zs{l\ge 1}x_\zs{l}
\E\,\int^{n}_\zs{0}\,I^{2}_\zs{t-}(\phi_\zs{j})\phi^{2}_\zs{l}(t)\d m_\zs{t}\\
&-
\sum_\zs{l\ge 1}x_\zs{l}
\E\,\int^{n}_\zs{0}\,
\left(\E\,I^{2}_\zs{t-}(\phi_\zs{j})\right)\phi^{2}_\zs{l}(t)\d m_\zs{t}
=0\,.
\end{align*}
Moreover, we note that
$$
\int^{n}_\zs{0}\,\wt{I}_\zs{t-}(\phi_\zs{j})\zeta_\zs{t-}(x)\d\xi_\zs{t}
=
\sum_\zs{l\ge 1}\,x_\zs{l}
\int^{n}_\zs{0}\,\wt{I}_\zs{t-}(\phi_\zs{j})
I_\zs{t-}(\phi_\zs{l})\,\phi_\zs{l}(t)
\,\d\xi_\zs{t}
$$
and
\begin{align*}
\int^{n}_\zs{0}\,\wt{I}_\zs{t-}(\phi_\zs{j})
I_\zs{t-}(\phi_\zs{l})\,\phi_\zs{l}(t)
\,\d\xi_\zs{t}
&=
\int^{n}_\zs{0}\,I^{2}_\zs{t-}(\phi_\zs{j})
I_\zs{t-}(\phi_\zs{l})\,\phi_\zs{l}(t)
\,\d\xi_\zs{t}\\[2mm]
&-
\int^{n}_\zs{0}\,\left(\E\,I^{2}_\zs{t}(\phi_\zs{j})\right)
I_\zs{t-}(\phi_\zs{l})\,\phi_\zs{l}(t)
\,\d\xi_\zs{t}\,.
\end{align*}
From Lemma~\ref{Le.sec:A.5}, it follows
$$
\E\,\int^{n}_\zs{0}\,\wt{I}_\zs{t-}(\phi_\zs{j})
I_\zs{t-}(\phi_\zs{l})\,\phi_\zs{l}(t)
\,\d\xi_\zs{t}=0
$$
and we come to \eqref{sec:Pr.16-1}.
Furthermore, by the Ito formula one obtains
\begin{align*}
\ov{I}^2_\zs{n}(x)&=2\int^n_\zs{0}\,\ov{I}_\zs{t-}(x)\d \ov{I}_\zs{t}(x)
+
4
\varrho^2_\zs{1}\int^n_\zs{0}\zeta^2_\zs{t}(x)\d t
+\varrho^{4}_\zs{2}\sum^{+\infty}_\zs{k=1}\,
\Phi^{2}_\zs{T_\zs{k}}(x)\,
Y^{4}_\zs{k}\Chi_\zs{\{T_\zs{k}\le n\}}\\
&+
4\varrho^{2}_\zs{2}
\sum^{+\infty}_\zs{k=1}\,
\zeta^{2}_\zs{T_\zs{k}-}(x)
Y^{2}_\zs{k}\,\Chi_\zs{\{T_\zs{k}\le n\}}
+4\varrho^{3}_\zs{2}
\sum^{+\infty}_\zs{k=1}\,
\Phi^{2}_\zs{T_\zs{k}}(x)
\zeta_\zs{T_\zs{k}-}(x)\,Y^{3}_\zs{k}
\Chi_\zs{\{T_\zs{k}\le n\}}
\,.
\end{align*}
By Lemma~\ref{Le.sec:A.3} one has $\E\,(\zeta_\zs{T_\zs{k}-}|T_\zs{k})=0$.
Therefore, taking into account  \eqref{sec:Pr.16-1}, we calculate
\begin{equation} \label{sec:Pr.17}
\E\ov{I}^2_\zs{n}(x)=
4 \varrho^2_\zs{1}
\E
\int^n_\zs{0}\zeta^2_\zs{t}(x)\d t
+
\varrho^4_\zs{2}\,\E\,Y^4_\zs{1}\, D_\zs{1,n}(x)
+
4\varrho^2_\zs{2} D_\zs{2,n}(x)
\,,
\end{equation}
where
$$
D_\zs{1,n}(x)=
\sum^\infty_\zs{k=1}
\E
\Phi^2_\zs{T_\zs{k}}(x)\Chi_\zs{\{T_\zs{k}\le n\}}
\quad\mbox{and}\quad
D_\zs{2,n}(x)=\sum^\infty_\zs{k=1}
\E
\zeta^2_\zs{T_\zs{k}-}(x)\Chi_\zs{\{T_\zs{k}\le n\}}\,.
$$
By applying Lemma~\ref{Le.sec:A.2},
one has
\begin{align}\nonumber
\E
\int^n_\zs{0}\zeta^2_\zs{t}(x)\d t&=\sum_\zs{i,j}x_\zs{i}x_\zs{j} \int^n_\zs{0}
\phi_\zs{i}(t)\phi_\zs{j}(t)
\E\,I_\zs{t}(\phi_\zs{i})I_\zs{t}(\phi_\zs{j}) \d t\\
\nonumber &=\sigma^* \sum_\zs{i,j}x_\zs{i}x_\zs{j} \int^n_\zs{0}
\phi_\zs{i}(t)\phi_\zs{j}(t) \left(
\int^t_\zs{0}\phi_\zs{i}(s)\phi_\zs{j}(s)\d s \right) \d t\\
\label{sec:Pr.18} &=\frac{\sigma^*}{2}
\sum_\zs{i,j}x_\zs{i}x_\zs{j} \left( \int^n_\zs{0}
\phi_\zs{i}(t)\phi_\zs{j}(t)\d t\right)^2 \le
n^2\frac{\sigma^*}{2} \,.
\end{align}
Further it is easy to check that
$$
D_\zs{1,n}=\lambda\int^n_\zs{0}\,\Phi^2_\zs{t}(x)\d t
=\lambda
\int^n_\zs{0} \left( \sum_\zs{j\ge 1}x_\zs{j}\phi^2_\zs{j}(t)
\right)^2\d t
\,.
$$
Therefore, taking into account that $\#(x)\le n$ and $|x|\le 1$,
we estimate $D_\zs{1,n}$ by applying the Causchy-Schwarts-Bounyakovskii inequality
\begin{equation}\label{sec:Pr.20}
D_\zs{1,n}\,
\le\, 4\lambda n \left( \sum_\zs{j\ge 1}x_\zs{j}
\right)^2
\le 4\lambda n \#(x)\le 4\lambda
n^2\,.
\end{equation}
Finally, we write down the process $\zeta_\zs{t}(x)$ as
$$
\zeta_\zs{t}(x)=\int^t_\zs{0}Q_\zs{x}(t,s)\d \xi_\zs{s} \quad
\mbox{with}\quad
  Q_\zs{x}(t,s)= \sum_\zs{j\ge
1}x_\zs{j}\phi_\zs{j}(s)\phi_\zs{j}(t)\,.
$$
By putting
$$
\wt{D}_\zs{2,n}= \E\, \sum^\infty_\zs{k=2}
\,\Chi_\zs{\{T_\zs{k}\le n\}}\, \sum^{k-1}_\zs{l=1}\,
Q^2_\zs{x}(T_\zs{k},T_\zs{l})
$$
and applying  Lemma~\ref{Le.sec:A.3} we obtain
\begin{align*}
D_\zs{2,n}&= \varrho^2_\zs{1} \sum^{\infty}_\zs{k=1}
\E\int^{T_\zs{k}}_\zs{0} Q^2_\zs{x}(T_\zs{k},s)\d s
\Chi_\zs{\{T_\zs{k}\le n\}} + \varrho^2_\zs{2}
\wt{D}_\zs{2,n}\\
&= \lambda\varrho^2_\zs{1} \int^n_\zs{0}\int^t_\zs{0}Q^2_\zs{x}(t,s)\d s
\d t + \varrho^2_\zs{2} \wt{D}_\zs{2,n} \,.
\end{align*}
Moreover, one can rewrite the second term in the last equality as
\begin{align*}
\wt{D}_\zs{2,n}&
=\sum^\infty_\zs{l=1}\E\,\Chi_\zs{\{T_\zs{l}\le n\}}
\sum^\infty_\zs{k=l+1}Q^2_\zs{x}(T_\zs{k},T_\zs{l})\Chi_\zs{\{T_\zs{k}\le n\}}\\
&=\lambda^2\int^n_\zs{0}\left(\int^{n-s}_\zs{0}Q^2_\zs{x}(s+z,s)\d z\right) \d s\\
&=\lambda^2\int^n_\zs{0}\left(\int^{t}_\zs{0}Q^2_\zs{x}(t,s)\d s\right) \d t
\,.
\end{align*}
Thus,
\begin{align}\nonumber
D_\zs{2,n}&\le (\lambda\varrho^2_\zs{1}+\lambda^2\varrho^2_\zs{2})
\int^n_\zs{0}\left(\int^{n}_\zs{0}Q^2_\zs{x}(t,s)\d s\right) \d t\\[2mm]
\label{sec:Pr.19}
&=
(\lambda\varrho^2_\zs{1}+\lambda^2\varrho^2_\zs{2})n^2
=\lambda\sigma^{*}n^{2}\,.
\end{align}
The equation \eqref{sec:Pr.17} and the inequalities
\eqref{sec:Pr.18}--\eqref{sec:Pr.20} imply the validity of
condition $\C_\zs{2})$ for the process \eqref{sec:In.3}.
Hence Theorem~\ref{Th.sec:2.2}.
\endproof

\subsection{Proof of Proposition~\ref{Pr.sec:Si.1}}\label{subsec:Pr.3}
Substituting \eqref{sec:Mo.4} in \eqref{sec:Si.1} yields
\begin{equation}\label{sec:Pr.21}
\wh{\sigma}_\zs{n}=\sum^n_\zs{j=l}\theta^2_\zs{j}
+\frac{2}{\sqrt{n}}\sum^n_\zs{j=l}\theta_\zs{j}\xi_\zs{j,n}
+\frac{1}{n}\sum^n_\zs{j=l}\xi^2_\zs{j,n}\,.
\end{equation}
Further, denoting
$$
x'_\zs{j}=\Chi_\zs{\{l\le j\le n\}}
\quad\mbox{and}\quad
x''_\zs{j}=\frac{1}{\sqrt{n}}\,\Chi_\zs{\{l\le j\le n\}}\,,
$$
we represent the last term in \eqref{sec:Pr.21} as
$$
\frac{1}{n}\sum^n_\zs{j=l}\xi^2_\zs{j,n}=
\frac{1}{n}B_\zs{1,n}(x')
+
\frac{1}{\sqrt{n}}
\,B_\zs{2,n}(x'')
+\frac{n-l+1}{n}\sigma\,,
$$
where the functions $B_\zs{1,n}(\cdot)$ and $B_\zs{2,n}(\cdot)$ are defined in
 conditions $\C_\zs{1})$ and $\C_\zs{2})$.
Combining these equations leads to the inequality
\begin{align*}
\E_\zs{S}|\wh{\sigma}_\zs{n}-\sigma|&\le \sum_\zs{j\ge l}\theta^2_\zs{j}
+
\frac{2}{\sqrt{n}}
\E_\zs{S}|\sum^n_\zs{j=l}\theta_\zs{j}\xi_\zs{j,n}|\\[2mm]
&+
\frac{1}{n}|B_\zs{1,n}(x')|
+
\frac{1}{\sqrt{n}}
\,\E\,|B_\zs{2,n}(x'')|
+\frac{l-1}{n}\sigma\,.
\end{align*}
By
Lemma~\ref{Le.sec:A.6} and
 conditions $\C_\zs{1})$, $\C_\zs{2})$, one gets
\begin{align*}
\E_\zs{S}|\wh{\sigma}_\zs{n}-\sigma|&\le \sum_\zs{j\ge l}\theta^2_\zs{j}
+
\frac{2}{\sqrt{n}}
\E_\zs{S}|\sum^n_\zs{j=l}\theta_\zs{j}\xi_\zs{j,n}|\\[2mm]
&+
\frac{c^{*}_\zs{1}(n)}{n}
+
\frac{c^{*}_\zs{2}(n)}{\sqrt{n}}
+\frac{\sigma}{\sqrt{n}}\,.
\end{align*}
In view of the inequality \eqref{sec:In.2}, the last term can be
estimated as
$$
\E_\zs{S}|\sum^n_\zs{j=l}\theta_\zs{j}\xi_\zs{j,n}|
\le
\sqrt{\sigma^*\sum^n_\zs{j=l}\theta^2_\zs{j}}\le
\sqrt{\sigma^*}|\dot{S}|_\zs{1}
\frac{2}{\sqrt{l}}\,.
$$
Hence Proposition~\ref{Pr.sec:Si.1}.
\endproof

\medskip

\renewcommand{\theequation}{A.\arabic{equation}}
\renewcommand{\thetheorem}{A.\arabic{theorem}}
\renewcommand{\thesubsection}{A.\arabic{subsection}}
\section{Appendix}\label{sec:A}
\setcounter{equation}{0}
\setcounter{theorem}{0}

\subsection{Property of the penalty term \eqref{sec:Mo.12}}\label{subsec:A.1}

\begin{lemma}\label{Le.sec:A.1}

Assume that the condition $\C_\zs{1})$ holds with $\sigma>0$.
Then for any $n\ge 1$ and $\gamma\in\Gamma$,
\begin{align*}
P_\zs{n}(\gamma)
\le \E_\zs{S}\,\Er_\zs{n}(\gamma)+
\frac{c^*_\zs{1}(n)}{n}\,.
\end{align*}
\end{lemma}
\proof By the definition of $\Er_\zs{n}(\gamma)$
 one has
\begin{align*}
\Er_\zs{n}(\gamma)=
\sum^\infty_\zs{j=1}
\left((\gamma(j)-1)\theta_\zs{j}+\gamma(j)
\frac{1}{\sqrt{n}}\xi_\zs{j,n}\right)^2\,.
\end{align*}
In view of the condition $\C_\zs{1})$ this leads to the desired
result
\begin{align*}
\E_\zs{S} \Er_\zs{n}(\gamma)\ge
\,\frac{1}{n}\sum^\infty_\zs{j=1}\gamma^2(j)\,\E\,\xi^2_\zs{j,n}
= P_\zs{n}(\gamma)-\frac{c^*_\zs{1}(n)}{n}\,.
\end{align*}
\endproof

\subsection{Properties of the process \eqref{sec:In.3}}\label{subsec:A.2}

\begin{lemma}\label{Le.sec:A.2}
Let $f$ and $g$ be any non-random functions from $\cL_\zs{2}[0,n]$ and
$(I_\zs{t}(f))_\zs{t\ge 0}$  be the process defined by
\eqref{sec:In.3}. Then, for any $0\le t\le n$,
$$
\E I_\zs{t}(f)I_\zs{t}(g)=\sigma^*\int^t_\zs{0}f(s)g(s)\d s
$$
where $\sigma^*=\varrho^2_\zs{1}+\lambda\varrho^2_\zs{2}$.
\end{lemma}
This Lemma is a direct consequence of Ito's formula as well as the
following result.
\begin{lemma}\label{Le.sec:A.3}
Let $Q$ be a bounded $[0,\infty)\times\Omega\to\bbr$ function
 measurable with respect to $\cB[0,+\infty)\bigotimes \cG_\zs{k}$,
where
\begin{equation}\label{sec:A.1}
\cG_\zs{k}=\sigma\{T_\zs{1},\ldots,T_\zs{k}\}
\quad\mbox{with some}\quad k\ge 2\,.
\end{equation}
Then
$$
\E\left( I_\zs{T_\zs{k-}}(Q)
|\cG_\zs{k}
\right)
=0
$$
and
$$
\E\left( I^2_\zs{T_\zs{k-}}(Q)
|\cG_\zs{k}
\right)
=\varrho^2_\zs{1}\int^{T_\zs{k}}_\zs{0}\,Q^2(s)\d s
+
\varrho^2_\zs{2}
\sum^{k-1}_\zs{l=1}Q^{2}(T_\zs{l})\,.
$$
\end{lemma}

Now we will study stochastic cadlag processes
$\eta=(\eta_\zs{t})_\zs{0\le t\le n}$ of the form
\begin{equation}\label{sec:A.2}
 \eta_\zs{t}=\sum^{\infty}_\zs{l=0}\,\upsilon_\zs{l}(t)\,
\Chi_\zs{\{T_\zs{l}\le t< T_\zs{l+1}\}}\,,
\end{equation}
where $\upsilon_\zs{0}(t)$ is a function measurable with respect
to $\sigma\{w_\zs{s}\,,\,s\le t\}$ and the coefficient
$\upsilon_\zs{l}(t)$, $l\ge 1$, is a function measurable with
respect to
$$
\sigma\{w_\zs{s},s\le
t\,,\,Y_\zs{1},\ldots,Y_\zs{l},T_\zs{1},\ldots,T_\zs{l}\}\,.
$$
Now we show the following result.
\begin{lemma}\label{Le.sec:A.4}
Let  $\eta=(\eta_\zs{t})_\zs{0\le t\le n}$ be a stochastic
non-negative process given by \eqref{sec:A.2}, such that
$$
\E\, \int^{n}_\zs{0}\,\eta_\zs{u}\,\d u\,<\,\infty\,.
$$
Then
$$
\E\,\int^{n}_\zs{0}\,\eta_\zs{u-}\,\d m_\zs{u}=0
$$
where the process $m=(m_\zs{t})$ is defined in \eqref{sec:Pr.16}.
 \end{lemma}
\proof Note that the stochastic integral, with respect to the
martingale \eqref{sec:Pr.16}, can be written as
\begin{align*}
\int^{n}_\zs{0}\,\eta_\zs{u-}\,\d m_\zs{u}&= \sum_\zs{0\le u\le
n}\,\eta_\zs{u-}\,(\Delta z_\zs{u})^{2} - \lambda
\int^{n}_\zs{0}\eta_\zs{u}\d u\\
& = \sum^{+\infty}_\zs{k=1}\,\eta_\zs{T_\zs{k}-}\,Y^{2}_\zs{k}\,
\Chi_\zs{\{T_\zs{k}\le n\}} - \lambda \int^{n}_\zs{0}\eta_\zs{u}\d
u \,.
\end{align*}
Therefore, taking into account the representation \eqref{sec:A.2},
we obtain
\begin{equation}\label{sec:A.3}
\int^{n}_\zs{0}\,\eta_\zs{u-}\,\d m_\zs{u}=\Upsilon_\zs{1}-\lambda\Upsilon_\zs{2}
\end{equation}
where
$$
\Upsilon_\zs{1}=\sum^{+\infty}_\zs{k=1}\,\upsilon_\zs{k-1}(T_\zs{k}-)\,Y^{2}_\zs{k}\,
\Chi_\zs{\{T_\zs{k}\le n\}}
\quad\mbox{and}\quad
\Upsilon_\zs{2}=\int^{n}_\zs{0}\eta_\zs{u}\d
u\,.
$$
Recalling that $\E Y^{2}_\zs{1}=1$ and
$\upsilon_\zs{k}\ge 0$, we calculate
$$
\E\Upsilon_\zs{1}=
\sum^{+\infty}_\zs{k=1}\,\E\,\upsilon_\zs{k-1}(T_\zs{k}-)\,
\Chi_\zs{\{T_\zs{k}\le n\}}\,.
$$
Moreover, the functions $(\upsilon_\zs{k})$ are cadlag processes, therefore
the Lebesgue measure of the set
$\{t\in\bbr_\zs{+}\,:\,\upsilon_\zs{k}(t-)\neq  \upsilon_\zs{k}(t)\}$
equals zero. Thus,
$$
\E\,\upsilon_\zs{k-1}(T_\zs{k}-)\, \Chi_\zs{\{T_\zs{k}\le n\}} =
\lambda\,\, \E\,\Chi_\zs{\{T_\zs{k-1}\le n\}}
\int^{n-T_\zs{k-1}}_\zs{0}\,\upsilon_\zs{k-1}(T_\zs{k-1}+u) \,
e^{-\lambda u}\,\d u\,.
$$
This implies
\begin{equation}\label{sec:A.4}
\E\Upsilon_\zs{1}=\lambda\sum^{+\infty}_\zs{l=0}
\E\,\Chi_\zs{\{T_\zs{l}\le n\}} \int^{n-T_\zs{l}}_\zs{0}
\,\upsilon_\zs{l}(T_\zs{l}+u) \, e^{-\lambda u}\,\d u\,.
\end{equation}
Similarly we obtain
\begin{align}\nonumber
\E\Upsilon_\zs{2}&=\sum^{+\infty}_\zs{l=0}\E\,\Chi_\zs{\{T_\zs{l}\le
n\}} \int^{n}_\zs{T_\zs{l}}\,\upsilon_\zs{l}(t)\, \Chi_\zs{\{t\le
T_\zs{l+1}\}} \d t\\ \label{sec:A.5}
 &= \sum^{+\infty}_\zs{l=0}
\E\,\Chi_\zs{\{T_\zs{l}\le n\}} \int^{n-T_\zs{l}}_\zs{0}
\,\upsilon_\zs{l}(T_\zs{l}+u) \, e^{-\lambda u}\,\d u\,.
\end{align}
Substituting \eqref{sec:A.4} and \eqref{sec:A.5} in \eqref{Le.sec:A.3}
implies the assertion of Lemma~\ref{Le.sec:A.4}.
\endproof

\begin{lemma}\label{Le.sec:A.5}
Assume that $\E\,Y^4_\zs{1}<\infty$. Then, for any measurable
bounded non-random functions $f$ and $g$, one has
$$
\E\,\int^{n}_\zs{0}\,I^{2}_\zs{t-}(f)I_\zs{t-}(g)\,g(t)\,\d\xi_\zs{t}
=0\,.
$$
\end{lemma}
\proof First we note that
$$
\E\,\int^{n}_\zs{0}\,I^{2}_\zs{t-}(f)I_\zs{t-}(g)\,g(t)\,\d
z_\zs{t}=\E\,\sum_\zs{j\ge
1}\,I^{2}_\zs{T_\zs{j}-}(f)I_\zs{T_\zs{j}-}(g)\,g(T_\zs{j})\,
\Chi_\zs{\{T_\zs{j}\le n\}}\,\E\, Y_\zs{j}\,
 =0\,.
$$
Therefore,  to prove this lemma one has to show that
\begin{equation}\label{sec:A.6}
\E\, \int^{n}_\zs{0}\,I^{2}_\zs{t}(f)I_\zs{t}(g)\,g(t)\,\d
w_\zs{t} \,=0 \,.
\end{equation}
To this end we represent the stochastic integral $I_\zs{t}(f)$ as
$$
I_\zs{t}(f)=\varrho_\zs{1}\,I^{w}_\zs{t}(f)\,+\,
\varrho_\zs{2}\,I^{z}_\zs{t}(f)
 \,,
$$
where
$$
I^{w}_\zs{t}(f)=\int^{t}_\zs{0}\,f_\zs{s}\,\d w_\zs{s}
 \quad\mbox{and}\quad
I^{z}_\zs{t}(f)= \int^{t}_\zs{0}\,f_\zs{s}\,\d z_\zs{s}\,.
$$
Note that
$$
\E\,|I^{z}_\zs{t}(f)|^{4}\,\le\,\,M^{4}
\,\E\,Y^{4}_\zs{1}\,\E\,N^{2}_\zs{n}
\,=\,M^{4}\,\E\,Y^{4}_\zs{1}(\lambda n+\lambda^{2}n^{2})
\,<\,\infty\,,
$$
where
$$
M=\sup_\zs{0\le t\le n}\left(|f(t)|+|g(t)|\right)\,.
$$
Therefore, taking into account that the processes
$(w_\zs{t})$ and
$(z_\zs{t})$ are independent, we get
$$
\E\,
\int^{n}_\zs{0}\,I^{4}_\zs{t}(f)\,(I^{w}_\zs{t}(g))^{2}\,g(t)\,\d t
 \,<\,\infty\,,
$$
i.e.
$$
\E\, \int^{n}_\zs{0}\,I^{2}_\zs{t}(f)I^{w}_\zs{t}(g)\,g(t)\,\d
w_\zs{t} \,=\,0\,.
$$
Similarly, we obtain
$$
\E\,
\int^{n}_\zs{0}\,(I^{w}_\zs{t}(f))^{2}I^{z}_\zs{t}(g)\,g(t)\,\d
w_\zs{t} \,=\,0
$$
and
$$
\E\, \int^{n}_\zs{0}\,I^{w}_\zs{t}(f)\,I^{z}_\zs{t}(f)
I^{z}_\zs{t}(g)\,g(t)\,\d w_\zs{t} \,=\,0\,.
$$
Therefore, to show \eqref{sec:A.6} one has to check that
\begin{equation}\label{sec:A.7}
\E\, \int^{n}_\zs{0}\,\eta_\zs{t}\,
\d w_\zs{t} \,=\,0 \,,
\end{equation}
where
$$
\eta_\zs{t}=\left(I^{z}_\zs{t}(f)\right)^{2}I^{z}_\zs{t}(g)\,g(t)\,.
$$
Taking into account that the processes $(\eta_\zs{t})$
and $(w_\zs{t})$ are independent, we get
$$
\E\, \left| \int^{n}_\zs{0}\,\eta_\zs{t}\d w_\zs{t} \right|
 \le\,
\E\,\sqrt{
\int^{n}_\zs{0}\,
\eta^{2}_\zs{t}\,\d t}
\,\le\,
\sqrt{n}\,\E\,\sup_\zs{0\le t\le n}\,|\eta_\zs{t}|
\,.
$$
Here,
the last term can be estimated as
$$
\E\,
\sup_\zs{0\le t\le n}
|\eta_\zs{t}|\le M^{4}\left(
\sum^{N_\zs{n}}_\zs{j=1}|Y_\zs{j}|
\right)^{3}\,\le\,M^{4}\,\E|Y_\zs{1}|^{3}\,
\E\,N^{3}_\zs{n}<\infty\,.
$$
Hence the stochastic integral $\int^{n}_\zs{0}\,\eta_\zs{t}\d w_\zs{t}$ is an integrable
random variable and
$$
\E\,\int^{n}_\zs{0}\,\eta_\zs{t}\,
\d w_\zs{t}=
\E\,\E\left(\int^{n}_\zs{0}\,\eta_\zs{t}\,
\d w_\zs{t}|\eta_\zs{t}\,,0\le t\le n\right)=0\,.
$$
Thus we obtain the equality \eqref{sec:A.7} which implies
\eqref{sec:A.6}. Hence Lemma~\ref{Le.sec:A.5}.
\endproof

\subsection{Property of the Fourier coefficients}\label{subsec:A.3}
\begin{lemma}\label{Le.sec:A.6}
Suppose that the  function $S$ in \eqref{sec:In.1} is
differentiable and satisfies the condition \eqref{sec:Si.2}. Then
the Fourier coefficients \eqref{sec:Mo.2} satisfy the inequality
$$
\sup_\zs{l\ge 2}\, l\,\sum^\infty_\zs{j=l}\theta^2_\zs{j}\, \le\,
4\,|\dot{S}|^2_\zs{1}\,.
$$
\end{lemma}
{\bf Proof}. In view of \eqref{sec:Mo.1}, one has
$$
\theta_\zs{2p}=-
\frac{1}{\sqrt{2}\pi p}\int^1_\zs{0}
\,\dot{S}(t)\sin(2\pi pt)\d t\,.
$$
and
\begin{align*}
\theta_\zs{2p+1}&=
\frac{1}{\sqrt{2}\pi p}\int^1_\zs{0}
\,\dot{S}(t)(\cos(2\pi pt)-1)\d t\\
&=
-
\frac{\sqrt{2}}{\pi p}\int^1_\zs{0}
\,\dot{S}(t)\sin^2(\pi pt)\d t\,,
\quad p\ge 1\,.
\end{align*}
From here, it follows that,
for any $j\ge 2$
$$
\theta^2_\zs{j}\,
\le\,
 \frac{2}{j^2}|\dot{S}|^2_\zs{1}\,.
$$
Taking into account that
$$
\sup_\zs{l\ge 2}
l\sum_\zs{j\ge l}\frac{1}{j^{2}}\,\le \,2\,,
$$
we arrive at the desired result.
\endproof

\end{document}